\documentclass[11pt,a4paper]{amsart}
\usepackage{multicol}
\usepackage{ifthen}
\usepackage{fancyhdr}
\usepackage{tabularx}
\usepackage{anysize}

\usepackage{latexsym}

\usepackage{amsmath,amssymb,amsfonts,amscd,amsthm}

\newtheorem{theorem}{Theorem}

\usepackage[matrix,arrow]{xy}

\xyoption{all}

\newcommand{\ro}{E(f_1,f_2) := {\Biggl\{ (x,\theta) \in M \times N^I \ \Big|\
\begin{array}{c} \theta(0) = f_1(x); \\ \theta(1) = f_2(x) \end{array}
\Biggr\}}}

\title{Coincidence theory in arbitrary  codimensions: \\ the minimizing problem}

\begin{document}
\vskip1cm

\author{Ulrich Koschorke }
\date{September 2004}

\maketitle

\vskip2cm
Let \ $f_1, f_2 : M \longrightarrow N$ \ be two (continuous) maps between smooth connected manifolds \ $M$ \ and \ $N$ \ without boundary, of strictly positive dimensions \ $m$ \ and \ $n$, resp., \ $M$ being compact. We are interested in making the coincidence locus
$$
C (f_1, f_2) := \{ x \in M \ | \ f_1 (x) = f_2 (x)\}
$$
as small (or simple in some sense) as possible after possibly deforming \ $f_1$ \ and \ $f_2$ \ by a homotopy.

\noindent{\bf Question.} \ How large is the {\bf m}inimum number of {\bf c}oincidence {\bf c}omponents
$$
MCC (f_1, f_2) := \min \{ \# \pi_0 (C (f'_1, f'_2)) \ | \ f'_1 \sim f_1, f'_2 \sim f_2\} \ ?
$$
In particular, when does this number vanish, i.e. when can \ $f_1$ \ and \ $f_2$ \ be deformed away from one another?

This is a very natural generalization of one of the central problems of classical fixed point theory (where \ $M = N$ \ and \ $f_2$ \ = \ identity map): \ determine  the minimum number of fixed points among all maps in a given homotopy class (see  \cite{Br} and \cite{BGZ}, proposition 1.5). Note, however, that in higher codimensions \ $m - n > 0$ \ the coincidence locus is generically a closed \ $(m -n)$-manifold so that it makes more sense to count {\it pathcomponents} \ rather than points. Also the methods of (first order, singular) (co)homology will no longer be strong enough to capture the subtle geometry of coincidence manifolds.

In this lecture I will use the language of normal bordism theory (and a nonstabilized version thereof) to define and study lower bounds \ $N (f_1, f_2)$ (and \ $N^\# (f_1, f_2))$ \  for \ $MCC (f_1, f_2)$.

After performing an approximation we may assume that the map \ $(f_1, f_2) : M \to N \times N$ \ is smooth and transverse to the diagonal \ $\Delta = \{ (y, y) \in N \times N \ | \ y \in N\}$. Then the coincidence locus
$$
C \ = \ C (f_1, f_2) \ = \ (f_1, f_2)^{- 1} (\Delta)
$$
is a closed smooth \ $(m - n)$-dimensional manifold, equipped with

\noindent i) \ \ maps
\begin{equation*}
\begin{xy}
    \POS(93,13.5)*+{\ro},
    \POS(35.5,2) \ar (52.5,11.5),
    \POS(41.5,9)*+{\tilde{g}}
    \POS(60,10) \ar (60,3),
    \POS(45,-3)*+{g = \operatorname{incl}}
    \POS(63.5,6.5)*+{\operatorname{pr}},
    \POS(32,0)*+{C}, \POS(36,0) \ar (56,0), \POS(60,0)*+{M},
\end{xy}
\end{equation*}
where \ $\widetilde g$ \ is the natural lifting which adds the constant path at \ $f_1 (x) = f_2 (x)$ \ to \ $g (x) = x \in C$; \ \ and

\noindent ii) \ a stable vector bundle isomorphism
$$
\overline g \ : \ TC \oplus g^* (f_1^* (TN)) \ \cong \ g^* (TM)
$$
deduced from the isomorphism
$$
\overline g^\# \ : \ \nu (C, M) \ \cong \ (f_1, f_2)^* (\nu (\Delta, N \times N)) \ \cong f_1^* (TN) \ | C
$$
of (nonstable) normal bundles.

The triple \ $(C, \widetilde g, \overline g)$ \ gives rise to a welldefined bordism class
$$
\widetilde\omega (f_1, f_2) \ := \ [C, \widetilde g, \overline g] \ \in \ \Omega_{m -n} (E (f_1, f_2); \ \widetilde\varphi)
$$
in the normal bordism group of such triples (here the virtual coefficient bundle is defined by
$$
\widetilde\varphi \ := \ pr^* (f^*_1 (TN) - TM) \ ;
$$
e.g.\ if \ $M$ \ and \ $N$ \ are stably parallelized, then \ $\widetilde\varphi$ \ is trivial and we are dealing with (stably) framed bordism).

Keeping track also of the fact that \  $C$ \ is a smooth submanifold of \ $M$ \ with (non-stabilized) normal bundle described by \ $\overline g^\#$, \ we obtain a sharper invariant
$$
\omega^\# (f_1, f_2) \ \in \ \Omega^\# (f_1, f_2)
$$
which, however, lies in general only in a suitable bordism  {\it set} \ (not group).

A crucial ingredient of both the \ $\widetilde\omega$- and the \ $\omega^\#$-invariant is the map \ $\widetilde g$. Indeed, the path space \ $E (f_1, f_2)$ \ has a very rich topology. Already its set \ $\pi_0 (E (f_1, f_2))$ \ of pathcomponents can be huge -- it corresponds bijectively to the so called Reidemeister set \ $R (f_1, f_2)$, a well-studied  set-theoretic quotient of the fundamental group \ $\pi_1 (N)$. This leads to a natural decomposition
$$
C (f_1, f_2) \ = \ \coprod_{A \in \pi_0 (E (f_1, f_2))} \ \widetilde g^{- 1} (A) \ \ .
$$
Let \ $N (f_1, f_2)$, and \ $N^\# (f_1, f_2)$, resp., denote the corresponding number of nontrivial contributions by the various pathcomponents \ $A$ \ of \ $E (f_1, f_2)$ \ to \ $\widetilde\omega (f_1, f_2)$ \ and \ $\omega^\# (f_1, f_2)$, resp.

\
\begin{theorem}
 \ \ {\rm(i)} \ \ The integers \ $N (f_1, f_2)$ \ and \ $N^\# (f_1, f_2)$ \ depend only on the homotopy classes of \ $f_1$ and $f_2$;

\noindent {\rm(ii)} \ \ $N (f_1, f_2) = N (f_2, f_1)$ \ and \ $N^{\#} (f_1, f_2) = N^\# (f_2, f_1)$;

\noindent {\rm(iii)} \ $0 \le N (f_1, f_2) \le N^\# (f_1, f_2) \le MCC (f_1, f_2) < \infty$;

\noindent {\rm(iv)} \ \ if \ $m = n$ \ then \ $N (f_1, f_2) = N^\# (f_1, f_2)$ \ coincides with the classical Nielsen number $($which has a standard definition at least if both \ $M$ \ and \ $N$ \ are orientable or if \ $f_2$ \ is the identity map$)$.
\end{theorem}

Recall the decisive progress made by \ {\it J.\ Nielsen} \ on the classical minimizing problem when he decomposed fixed point sets into equivalence classes. In our interpretation this is just the decomposition of a \ $0$-dimensional bordism class according to the pathcomponents of its target space. In higher (co)dimensions \ $m -n$ \ the map \ $\widetilde g$ \ into \ $E (f_1, f_2)$ \ and the twisted framing  \ $\overline g^\#$ \ contain much more information. E.g.\ if \ $M = S^m$ \ and \ $n \ge 2$ \ then \ $\Omega^\# (f_1, f_2)$ \ can be identified  with the homotopy group \ $\pi_m (S^n \wedge \Omega (N)^+)$, and \ $\omega^\# (f_1, f_2)$ \ is closely related to a Hopf-Ganea invariant. This allows us to reduce many aspects of our problem to questions in standard homotopy theory.

Details of definitions, proofs, and applications will be given elsewhere (compare e.g. \cite{K3} and \cite{K2}). Here we present just one sample result.

\begin{theorem} \ Let \ $N$ \ be an odd-dimensional spherical  space form $($i.e.\ the quotient of \ $S^n$ \ by a free action of a finite group$)$. Then we have for all \ $f_1, f_2 : S^m \to N$:
$$
MCC (f_1, f_2) = N^\# (f_1, f_2) = \left\{ \begin{array}{ll}
0 & \text{if} \ f_1 \sim f_2 \ \text{or} \ m < n; \\
\#\pi_1 (N) & \text{if} \ f_1 \not\sim f_2 \ \text{and} \ m > 1; \\
|d^0 (f_1) - d^0 (f_2)| & \text{if} \ m = n = 1 . \end{array} \right.
$$
$($Here \ $d^0 (f_i) \in  \mathbb Z$ \ denotes the usual degree$)$.
\end{theorem}

Finally note that our approach applies also to over- and undercrossings of link maps into a manifold of the form \ $N \times \mathbb R$. This yields unlinking obstructions which often settle unlinking questions and which, in addition, turn out to distinguish a great number of different link homotopy classes. In certain cases they even allow a complete link homotopy classification. Moreover, our approach leads also to the notion of Nielsen numbers of link maps (cf.\ \cite{K4}).

\providecommand{\bysame}{\leavevmode\hbox to3em{\hrulefill}\thinspace}
\providecommand{\MR}{\relax\ifhmode\unskip\space\fi MR }
\providecommand{\MRhref}[2]{%
  \href{http://www.ams.org/mathscinet-getitem?mr=#1}{#2}
}
\providecommand{\href}[2]{#2}

\vspace*{0.5cm}
\noindent
The last preprints are available at \newline http://www.math.uni-siegen.de/topology/publications.html.


\vfill

%



\end{document}